\documentclass{article}

\def\d{ \delta }

\def\1ox{{ \Omega^1_{\scriptstyle{X}} }}
\def\2ox{{ \Omega^2_{\scriptstyle{X}} }}

\def\ok1{{ \Omega^1_K }}
\def\ok2{{ \Omega^2_K }}
\def\Om{{ \Omega }}
\def\om{{ \omega  }}
\def\O{{ {\mathcal O} }}

\def\ra{{ \rightarrow }}

\def\g{{ \gamma }}

\def\8{{ {\infty } }}

\def\deg{{ \mbox{deg} }}
\def\^{{ ^{\wedge} }}

\newtheorem{thm}{Theorem}

\usepackage{amsmath}
\usepackage{amsfonts}
\usepackage{amscd}
\usepackage{amssymb}

\title{A note on Szpiro's inequality for curves of higher genus
}
\author{ Minhyong Kim}

\begin{document}
\maketitle

Let $f:X \ra B$ be a semi-stable family of curves of genus $g$
over
a smooth projective curve $B$ of genus $\g$
over the complex numbers.
Assume that the family is relatively minimal in that
there are no $(-1)$ curves contained in the fibers and that it is
non-isotrivial.
Let $\d$ be the number of singular points in the
fibers and $s$ the number of singular fibers, i.e.,
the {\em singular values} of the map.
A. Beauville \cite{Be} has pointed out the following inequality
$$\d < (4g+2)(2\g-2+s)$$
which generalizes the `conductor-discriminant inequality'
for elliptic curves whose significance was pointed out by Szpiro \cite{Sz}.

This inequality is a quick consequence of the Bogomolov-Miyaoka-Yau
inequality and an inequality of Xiao \cite{Xi}
$$\d \leq (8+4/g)\deg(f_*(\om_{X/B}))$$
where $\om_{X/B}$ is the relative dualizing sheaf of the fibration.

The purpose of this note is to point out the following variant:
\begin{thm}
Let $g_0$ be the dimension of the constant part of the
Jacobian fibration of $X/B$. Then
$$\d \leq (g-g_0)(4+2/g)(2\g-2+s)$$
\end{thm}
Note that when $g_0=0$, this reduces to a slightly
weaker version of Beauville's inequality (in that it is not
strict). Even though this is just a minor variation
on Beauville's result we judged it worth stating since
the B-M-Y inequality is not used in the proof.
On the one hand, this makes it more likely to be extended to
positive characteristic. On the other, the ideas involved
may be more likely
to fit into the current interest in finding an`arithmetic'
Kodaira-Spencer map.
That is, the proof is similar to Beauville's
except we replace the B-M-Y inequality by the following inequality \cite{K}
$$\deg(f_*(\om_{X/B}))\leq ((g-g_0)/2)(2\g-2+s)$$
We recall briefly the proof of this inequality, referring to op. cit.
for details. There it is stated for abelian varieties, so it is
useful to quickly point out  the direct (and easier) argument for curves.:
We have the logarithmic Kodaira-Spencer map of
the fibration
$$KS: f_*(\om_{X/B}) \ra R^1f_*(\O_X) \otimes \Om^1_B(S)$$
where $S\subset B$ denotes the divisor of singular values.
By duality, one has $$R^1f_*(\O_X)\simeq (f_*(\om_{X/B}))^*$$
On the other hand, if $K$ denotes the kernel of $KS$,
a simple computation using the constancy of the  duality
pairing w.r.t. the Gauss-Manin connection shows
that \cite{K}
$$Im(KS) \simeq [f_*(\om_{X/B})/K]^*$$
So we have a full-rank map
$$f_*(\om_{X/B})/K \ra [f_*(\om_{X/B})/K]^* \Om^1_B(S)$$
and hence, taking determinants,
a non-zero map
$$[\Lambda^{(g-g_0)} (f_*(\om_{X/B})/K )]^{\otimes 2} \ra 
(\Om^1_B(S))^{\otimes(g-g_0)}$$
>From this, we get the inequality
$$\deg [f_*(\om_{X/B})/K ] \leq [(g-g_0)/2 ](2\g-2+s)$$
However, $\deg(K)=0$. One proves this by showing that
this kernel is preserved by the Gauss-Manin connection.
This implies that its top exterior power  also
carries a connection. But since the Gauss-Manin connection
is unipotent near $S$, we get that the connection on the top exterior
power is actually regular on all of $B$. Therefore,
the determinant of $K$ must have vanishing first Chern class.

These considerations give us the desired inequality for
$\deg [f_*(\om_{X/B})]$.
By combining this with Xiaos' inequality, we get
$$\d \leq (8+4/g)\deg [f_*(\om_{X/B})] \leq (8+4/g)[(g-g_0)/2](2\g-2+s)$$
which is the stated inequality.

\medskip
{\bf Acknowledgement:} I am grateful to Professor Beauville
for a
quick
reponse to an e-mail inquiry that prompted this note.

{\footnotesize DEPARTMENT OF MATHEMATICS, UNIVERSITY OF ARIZONA, TUCSON, AZ 85721, EMAIL: kim@math.arizona.edu and
KOREA INSTITUTE FOR ADVANCED STUDY, 207-43 CHEONGRYANGRI-DONG, DONGDAEMUN-GU,
SEOUL, KOREA 130-012}


\begin{thebibliography}{10}
\bibitem{Be} Beauville, Arnaud The Szipro inequality for higher genus fibrations.
Preprint (2001). math.AG/0109080
\bibitem{K} Kim, Minhyong ABC inequalities for some moduli spaces of
log general type, Math. Res. Letters 5 (1998), no.4, 517-522.
 
\bibitem{Sz} Szpiro, Lucien Discriminant and conducteur des courbes
elliptques. {\em S\'eminaire sur les pinceaux des courbes elliptiques}
(Paris 1988), Ast\'erisque 183 (1990), 7-18.
\bibitem{Xi} Xiao, Gang Fibered algebraic surfaces with low slope.
Math. Ann. 276 (1987), 449-466.

\end{thebibliography}
\end{document}